\documentclass[12pt,reqno]{amsart} 
\usepackage{amscd,amsmath,amsfonts,amsthm,epsfig, graphics, latexsym, mathrsfs}
\usepackage{amssymb}
\usepackage[unicode]{hyperref}

\usepackage[alphabetic,initials]{amsrefs}

\theoremstyle{plain}

\newtheorem{thm}{Theorem}[section]
\newtheorem{prop}[thm]{Proposition}
\newtheorem{lemma}[thm]{Lemma}
\newtheorem{cor}[thm]{Corollary}

\newtheorem{ques}[thm]{Question}

\newtheorem{add}[thm]{Addendum}

\theoremstyle{definition}
\newcommand{\appsection}[1]{\let\oldthesection\thesection
\renewcommand{\thesection}{Appendix \oldthesection}
\section{#1}\let\thesection\oldthesection}
\newtheorem{defn}[thm]{Definition}

\newtheorem{say}[thm]{}
\newtheorem{example}[thm]{Example}

\theoremstyle{remark}

\newtheorem{rmk}[thm]{Remark}
\newtheorem{rmks}[thm]{Remarks}

\newcommand{\bi}{\begin{itemize}}  
\newcommand{\ei}{\end{itemize}}
\newcommand{\bp}{\begin{proof}}
\newcommand{\ep}{\end{proof}}

\def\AA{{\mathbb{A}}}
\def\PP{{\mathbb{P}}}
\def\RR{{\mathbb{R}}}
\def\ZZ{{\mathbb{Z}}}
\def\QQ{{\mathbb{Q}}}
\def\CC{{\mathbb{C}}}

\def\GG{{\mathbb{G}}}

\def\kk{{\mathbf{k}}}

\def\gp{\mathfrak{p}}

\DeclareMathOperator{\Aut}{Aut}
\DeclareMathOperator{\Bl}{Bl}
\DeclareMathOperator{\Cl}{Cl}

\DeclareMathOperator{\Cox}{Cox}
\DeclareMathOperator{\Eff}{\overline{Eff}}
\DeclareMathOperator{\Hom}{Hom}
\DeclareMathOperator{\HH}{H}
\DeclareMathOperator{\Mov}{Mov}
\DeclareMathOperator{\NN}{N}
\DeclareMathOperator{\NE}{\overline{NE}}
\DeclareMathOperator{\Nef}{Nef}
\DeclareMathOperator{\PGL}{PGL}
\DeclareMathOperator{\Pic}{Pic}
\DeclareMathOperator{\Proj}{Proj}
\DeclareMathOperator{\R}{R}

\DeclareMathOperator{\Spec}{Spec}
\DeclareMathOperator{\WDiv}{WDiv}

\DeclareMathOperator{\ch}{char}

\DeclareMathOperator{\rk}{rk}

\newcommand{\cC}{\mathcal C}

\newcommand{\cO}{\mathcal O}

\newcommand{\cW}{\mathcal W}

\newcommand{\al}{\alpha}
\newcommand{\be}{\beta}
\newcommand{\de}{\delta}

\newcommand{\Ga}{\Gamma}
\newcommand{\De}{\Delta}
\newcommand{\La}{\Lambda}
\newcommand{\Si}{\Sigma}

\newcommand{\ra}{\rightarrow}

\newcommand{\dra}{\dashrightarrow}

\newcommand{\M}{\overline{M}}
\newcommand{\LM}{\overline{LM}}

\newcounter{et}[section]

\begin{document}
\bibliographystyle{amsplain}
\title{Mori Dream Spaces and Blow-ups}

\author{Ana-Maria Castravet}

\address{Ana-Maria Castravet: \sf Department of Mathematics, Northeastern University, 360 Huntington Avenue, Boston, MA 02115} 
\email{noni@alum.mit.edu}


\begin{abstract}
The goal of the present article is to survey  the general theory of Mori Dream Spaces, with special regards to the question: When is the blow-up of toric variety at a general point a Mori Dream Space? We translate the question for toric surfaces of Picard number one into an interpolation problem involving points in the projective plane. An instance of such an interpolation problem is the Gonzalez-Karu theorem that gives new examples of weighted projective planes whose blow-up at a general point is not a Mori Dream Space. 

\end{abstract}

\maketitle

\section{Introduction}

Mori Dream Spaces were introduced in \cite{HK} as a natural Mori theoretic generalization of toric varieties. As the name suggests, their main feature is that the Minimal Model Program (MMP) can be run for any divisor (not just the canonical divisor class). In particular, as for toric varieties, one only has to look into the combinatorics of the various birational geometry cones to achieve the desired MMP steps. 

As being a Mori Dream Space is equivalent to all (multi-)section rings being finitely generated, it is not surprising that non-trivial examples may be hard to find. It was not until the major advances in the MMP, that Hu and Keel's original conjecture that varieties of Fano type are Mori Dream Spaces was proved  \cite{BCHM}.  Although there are many examples outside of the Fano-type range, these often have an ad-hoc flavor. Certain positivity properties of the anticanonical divisor (such as being of Fano type or Calabi-Yau type)  of a Mori Dream Space are reflected in the multi-section rings \cite{Okawa}, \cite{GOST}, but no clear picture emerges in general.
More often than not, the usual operations of blowing up, taking projective bundles, crepant resolutions, hyperplane sections, when applied to Mori Dream Spaces, do not lead to Mori Dream Spaces. 

\medskip 

Our current goal is to pay special attention to blow-ups of Mori Dream Spaces, in particular, blow-ups at a single (general) point. More specifically, the following is a question asked by Jenia Tevelev:
\begin{ques}\label{Jenia}
Let $X$ be a projective $\QQ$-factorial toric variety over an algebraically closed field $\kk$. When is the blow-up $\Bl_p X$ of $X$ at a general point $p$ not a MDS?
\end{ques}

Using the action of the open torus $T=(\kk^*)^n$, we may assume the point $p$ is the identity $e$ of $T$.  Currently, the only known examples of $X$ toric such that $\Bl_e X$ is not a MDS fall into the following categories:

\bi
\item[(I) ] Certain (singular) toric projective surfaces with Picard number one;
\item[(II) ] Certain toric varieties  for which there exists a small modification that admits a surjective morphism  into one of the toric surfaces in (I). (Note that small modifications and images of Mori Dream Spaces are Mori Dream Spaces \cite{HK}, \cite{Okawa}).
\ei

All known examples are in characteristic zero, since the only examples of surfaces in (I) are in characteristic zero.  Eventually, blowing up (very) general points\footnote{Recall that the blow-up of a toric variety along a torus invariant stratum is a toric  variety.} on a toric variety leads to non Mori Dream Spaces: for example, the blow-up of $\PP^2$ at $r$ very general points is toric if and only if $r\leq 3$ and a Mori Dream Space if and only if $r\leq 8$.

\medskip 

A good portion of the examples in (I) are weighted projective planes $\PP(a,b,c)$ for a certain choice of weights $(a,b,c)$. Until \cite{CT4}, \cite{GK}, the only known examples of varieties as in Question \ref{Jenia} were of this type \cite{GNW}. The question whether $\Bl_e \PP(a,b,c)$ is a Mori Dream Space is equivalent to the symbolic Rees algebra of a so-called monomial prime ideal being Noetherian, and as such, it has a long history. Major progress was recently achieved by Gonzalez and Karu \cite{GK} by using methods of toric geometry. However, the main question remains open:

\begin{ques}\label{classify}
For which triples $(a,b,c)$  the blow-up $\Bl_e\PP(a,b,c)$ of $\PP(a,b,c)$ at the identity point $e$ is not a MDS?
\end{ques}

With the exception of $(a,b,c)=(1,1,1)$, in all examples where the Mori Dream Space-ness of $\Bl_e\PP(a,b,c)$ is understood (one way or another), it happens that $\Bl_e\PP(a,b,c)$ contains a negative curve $C$, different than the exceptional divisor $E$ above the point $e$. In positive characteristic, the existence of the negative curve $C$ implies that  $\Bl_e\PP(a,b,c)$ is a  Mori Dream Space by Artin's contractability theorem \cite{Artin}. No triples 
$(a,b,c)\neq(1,1,1)$ are known for which $\Bl_e\PP(a,b,c)$ contains no negative curve (other than $E$). If such an example exists (in any characteristic), it would imply the Nagata conjecture on linear systems on blow-ups of $\PP_{\CC}^2$ at $abc$ points \cite{CutkoskyKurano}.  If $\sqrt{abc}\notin\ZZ$, such an example would have many important consequences: new cases of the Nagata conjecture, examples of irrational Seshadri constants, new examples when $\Bl_e\PP(a,b,c)$ is not a Mori Dream Space, etc.

\

The goal of the present article is two-fold. First, to survey some of the general theory of Mori Dream Spaces, along with known results and open problems related to Question \ref{Jenia}. Second, to use the toric geometry methods of Gonzalez and Karu in order to translate Question \ref{classify} (and more generally, Question \ref{Jenia} in the case of surfaces of Picard number one) into an interpolation problem involving points in the (usual) projective plane $\PP^2$ (this translation is likely not new to the experts). As an illustration of this approach, we reprove (or rather, present a shortcut in the proof of) the main theorem in \cite{GK} (Thm. \ref{GK}). The advantages are that the interpolation problem is really equivalent to the original question, and there are further potential applications towards  
Question \ref{Jenia} and  Question \ref{classify}. For example, both of the following questions can be reformulated into interpolation problems: (a) whether $\Bl_e\PP(a,b,c)$ is a Mori Dream Space when in the presence of a negative curve, or (b)
whether  $\Bl_e\PP(a,b,c)$ has any negative curves at all. The drawback is that the interpolation problem seems to be almost equally difficult. 

By interpolation, we simply mean to separate points lying in the lattice points of a plane polytope (so in a grid!)  by curves of an appropriate degree. For example, to prove that $\Bl_e\PP(9,10,13)$ has no negative curve (other than $E$), it suffices to answer affirmatively:
 
\begin{ques}\label{interpolation_9,10,13}
Let $\De$ be the polytope in $\RR^2$ with vertices $(0,0)$, $(10,40)$, $(36,27)$. For every $q\geq 1$, let 
$$m_q=\lfloor q\sqrt{1170}\rfloor+1.$$
Is it true that for every $q\geq1$ and any point $(i,j)\in q\De\cap\ZZ^2$, there exists a curve $C\subset\RR^2$ of degree $m_q$ passing through all the points 
$(i',j')\neq(i,j)$ in $q\De\cap\ZZ^2$, but not $(i,j)$?
\end{ques}

\medskip

\noindent {\bf Structure of paper.} The first three sections present a general survey on Mori Dream Spaces: Section \ref{basics}
reviews the basic definitions and properties, Section \ref{examples} presents several key examples, while Section \ref{structure} gives an overview of the ``structure theory". The last four sections focus on blow-ups at a general point. Section \ref{Picard1} discusses generalities on blow-ups of 
(not necessarily toric) surfaces of Picard number one, while Section \ref{section wpp} presents the special case of weighted projective planes. Section \ref{higher} discusses blow-ups of higher dimensional toric varieties, with Losev-Manin spaces playing a central role. Finally (the linear algebra heavy) Section \ref{toricPicard1} translates Question \ref{Jenia} in the case of surfaces of Picard number one, into an interpolation problem and proves Thm. \ref{GK} as an application.

\medskip

\noindent {\bf Conventions and Notations.} Unless otherwise specified, we work over an algebraically closed field $\kk$ of arbitrary characteristic. 
For an abelian group $\Ga$ and a field $K$, we denote $\Ga_K$ the $K$-vector space $\Ga\otimes_{\ZZ} K$. 

\medskip

\noindent {\bf Acknowledgements.} I am grateful to Jenia Tevelev who pointed out Question \ref{Jenia} and the surrounding circle of ideas. 
I thank Shinosuke Okawa for his questions and comments, Jos\'e Gonzalez and Antonio Laface for useful discussions, and the 
anonymous referees for several useful comments. This work is partially supported by NSF grant DMS-1529735. I thank Institut de 
Math\'ematiques de Toulouse for its hospitality during the writing of this paper.


\section{Mori Dream Spaces}\label{basics}

Mori Dream Spaces are intrinsically related to Hilbert's $14$'th problem. Many of the results on finite generation of multi-section rings go back to Zariski and Nagata (see \cite{Mumford}). For a survey of Mori Dream Spaces from the invariant theory perspective, see \cite{McK_survey}.
In what follows, we briefly recall the definitions and basic properties from \cite{HK}. We found \cite{Okawa} to be a useful additional reference.

 \medskip
 
Let $X$ be a projective variety over $\kk$. We denote by $\NN^1(X)$ the group of Cartier divisors modulo numerical equivalence\footnote{$\NN^1(X)$ is a finitely generated abelian group.}. The cone generated by nef divisors in $\NN^1(X)_{\RR}$ is denoted $\Nef(X)$. Similarly, the closure of the cone of effective divisors (resp., movable divisors) is denoted $\Eff(X)$ (resp., $\Mov(X)$). Recall that an effective divisor is called movable if its base locus has codimension at least $2$. Similarly, if $\NN_1(X)$ is the group of $1$-cycles modulo numerical equivalence\footnote{The dual of $\NN^1(X)$ under the intersection pairing.}, the Mori cone $\NE(X)$ is the closure in $\NN_1(X)_{\RR}$ of the cone of effective $1$-cycles. 

\smallskip

The closure operations in the definition of $\Eff(X)$, $\Mov(X)$ and $\NE(X)$ are not necessary for Mori Dream Spaces (see Prop. \ref{effective} below). 
A \emph{small $\QQ$-factorial modification} (SQM for short) of a normal projective variety $X$ is a small (i.e., isomorphic in codiemsnion one) birational map $X\dra Y$ to another normal, $\QQ$-factorial projective variety $Y$.

\begin{defn}\label{mds}
A normal projective variety $X$ is called a \emph{Mori Dream Space} (MDS for short) if the following conditions are satisfied:
\bi
\item[(1) ] $X$ is $\QQ$-factorial, $\Pic(X)$ is finitely generated, with
$$\Pic(X)_{\QQ}\cong\NN^1(X)_{\QQ};$$
\item[(2) ] $\Nef(X)$ is generated by finitely many semiample divisors; 
\item[(3) ] There are finitely many SQMs $f_i: X\dra X_i$ such that each $X_i$ satisfies (1) and (2), and $\Mov(X)$ is the union of $f_i^*\Nef(X_i)$\footnote{If $f: X\dra Y$ is birational map, the pull back $f^*D$ of a Cartier divisor $D$ from $Y$ is defined as $p_*(q^*D)$, where $p: W\ra X$, $q: W\ra Y$ are given by a common resolution. If $f$ is small, $f^*D$ is simply the push forward $f^{-1}_*(D)$ via the inverse map $f^{-1}$.}.
\ei
\end{defn}

\begin{rmks}
(a) If $\kk$ is not the algebraic closure of a finite field, the condition that $\Pic(X)$ is finitely generated is equivalent to the condition $\Pic(X)_{\QQ}\cong\NN^1(X)_{\QQ}$, but not otherwise (see \cite[Rmk. 2.4]{Okawa})\footnote{In the original definition in \cite{HK}, only the condition 
$\Pic(X)_{\QQ}\cong\NN^1(X)_{\QQ}$ appears, but as explained in \cite{Okawa}, adding both conditions seems more natural.}. 

(b) Semiampleness and polyhedrality in conditions (2) and (3) are key, guaranteeing that all the MMP steps are reduced to combinatorics (finding the divisor class with the desired numerical properties).
\end{rmks}

A birational map $f: X\dra Y$ between normal projective varieties is called contracting if the inverse map $f^{-1}$ does not contract any divisors. If $E_1,\ldots E_k$ are the prime divisors contracted by $f$, then $E_1,\ldots E_k$ are linearly independent in $\NN^1(X)_{\RR}$ and each $E_i$ spans an extremal ray of $\Eff(X)$. The effective cone of a MDS also has a decomposition into rational polyhedral cones:

\begin{prop}\cite[Prop. 1.11 (2)]{HK}\label{effective}
Let $X$ be a MDS. There are finitely many birational contractions $g_i: X \dra Y_i$, with $Y_i$ a MDS, such that
$$\Eff(X)=\bigcup_i\cC_i,$$
$$\cC_i=g_i^*\Nef(Y_i)+\RR_{\geq0}\{E_1,\ldots,E_k\},$$
where $E_1,\ldots, E_k$ are the prime divisors contracted by $g_i$. 
\end{prop}

The cones $\cC_i$ are called the \emph{Mori chambers} of $X$.
Prop. \ref{effective} is best interpreted as an instance of Zariski decomposition: for each effective $\QQ$-Cartier divisor $D$, there exists a birational contraction $g: X\dra Y$ (factoring through an SQM and a birational morphism $X\dra X'\ra Y$) and $\QQ$-divisors $P$ and $N$, such that $P$ is nef on $X'$, $N$ is an effective divisor contracted by $g$ and for $m>0$ sufficiently large and divisible, the multiplication map given by the canonical section $x^m_N$ 
$$\HH^0(X,\cO(mP))\ra\HH^0(X,\cO(mD))$$
is an isomorphism. To see this, simply take 
$$P=g^*g_*(D),\quad N=D-P.$$

\begin{rmks}\label{rmk2}\label{rmks1}
(a) If $X$ is a MDS, all birational contractions $X\dra Y$ with $\QQ$-factorial $Y$, are the ones that appear in Prop. \ref{effective}. 
In particular, any such $Y$ is a MDS.

(b) The SQMs in Def. \ref{mds} are the only SQMs of $X$.  In particular, any SQM of a MDS is itself a MDS. 
\end{rmks}

\begin{defn} Let $X$ be a normal variety. For a semigroup $\Ga\subset\WDiv(X)$\footnote{$\WDiv(X)$ is the group freely generated by prime Weil divisors in $X$.} of Weil divisors on $X$, we define the multi-section ring $\R(X,\Ga)$ as the $\Ga$-graded ring:
$$\R(X,\Ga)=\bigoplus_{D\in\Ga}\HH^0(X,\cO(D))$$
with the multiplication induced by the product of rational functions. When $\Ga$ is a group such that the class map $\Ga_{\QQ}\ra\Cl(X)_{\QQ}$ is an isomorphism, we call $\R(X,\Ga)$ a \emph{Cox ring} of $X$ and denote this by $\Cox(X)$\footnote{The greater generality of working with Weil divisors rather than than Cartier divisors will be essential in Section \ref{section wpp}.}. 
\end{defn}

The definition of $\Cox(X)$ depends on the choice of $\Ga$, but basic properties, such as finite generation as a $\kk$-algebra, do not.  Note that if $\Ga'\subset\Ga$ is finite index subgroup, then $\R(X,\Ga)$ is an integral extension of $\R(X,\Ga')$.  For more details on Cox rings see \cite{ADHL}, \cite{LafaceVelasco}. 

\medskip

Mori Dream Spaces can be algebraically characterized as follows:
\begin{thm}\cite[Prop. 2.9]{HK}\label{cox} 
Let $X$ be a projective normal variety  satisfying condition (1) in Def. \ref{mds}. Then $X$ is a MDS if and only if $\Cox(X)$ is a finitely generated $\kk$-algebra. 
\end{thm}

\bp[Sketch of Proof]
If $\Cox(X)$ is finitely generated, let $V$ be the affine variety $\Spec(\Cox(X))$. Since $\Cox(X)$ is
graded by a lattice $\Ga\subset\WDiv(X)$, the algebraic torus $T=\Hom(\Ga,\GG_m)$ naturally acts on the affine variety $V$.
Let $\chi\in\Ga$ be a character of $T$ which corresponds to an ample divisor in $\Ga$. Then $X$ is $V//_{\chi}T$, the GIT quotient 
constructed with respect to the trivial line bundle on V endowed with a $T$-linearization by $\chi$. Similarly, all small modifications of $X$ can be obtained as GIT quotients $V//_{\chi}T$, for different classes $\chi$ in $\Ga$ (thus the Mori chamber decomposition is an instance of variation of GIT). 
The ``only if" implication follows from the more general Lemma \ref{multisection}. 
\ep

\begin{lemma}\label{multisection} 
Let $X$ be a MDS and let $\Ga$ be a finitely generated group of Weil divisors. Then $\R(X,\Ga)$ is a finitely generated $\kk$-algebra. 
\end{lemma}

\bp We follow the proof in \cite[Lemma 2.20]{Okawa}. 
The key facts used are (i)  $\R(X,\Ga)$ is finitely generated if $\Ga$ is generated by finitely many semiample divisors (\cite[Prop. 2.8]{HK}); (ii) Zariski decomposition as in Prop. \ref{effective}.
When $\R(X,\Ga)$ is a Cox ring, this is immediate: as $\Nef(X)$ is a full cone inside $\NN^1(X)_{\RR}$, if $\Ga$ is generated by $\QQ$-divisors that are generators of $\Nef(X)$ (hence, $\Ga_{\QQ}\cong\Cl(X)_{\QQ}$), the result follows by (i). 

For the general case, without loss of generality, we may replace  $\Ga$ with a subgroup of finite index. In particular, we may assume that $\Ga$ has no torsion. For a Mori chamber $\cC$, denote $\Ga_{\cC}=\Ga\cap\cC$ (a semigroup). 
As there are finitely many Mori chambers and the support of $\R(X,\Ga)$ is the union of $\Ga_{\cC}$, it is enough to prove that $\R(X,\Ga_{\cC})$ is finitely generated.  We may assume that there is $g:X\ra Y$ birational morphism, with $$\cC=g^*\Nef(Y)+\RR_{\geq0}\{E_1,\ldots,E_k\},$$
where $E_1,\ldots,E_k$ are the prime divisors contracted by $g$. Note that since $\cC$ is a rational polyhedral cone,  $\Ga_{\cC}$ is a finitely generated semigroup. For a set of generators $D_1,\ldots, D_r$ we consider Zariski decompositions as in Prop. \ref{effective}:
$D_i=P_i+N_i$, with  $\QQ$-divisors  $P_i$ in $g^*\Nef(X)$ and $N_i$ effective and supported on $E_1,\ldots, E_k$. Up to replacing each $D_i$ with a multiple, we may assume $P_i$ and $N_i$ are $\ZZ$-divisors. Then  $\R(X,\Ga_{\cC})$ is isomorphic to an algebra over  $\R(Y,P_1,\ldots, P_r)$ generated by the canonical sections $x_{N_1},\ldots,x_{N_k}$. By (i), it follows that $\R(X,\Ga_{\cC})$ is finitely generated. 

\ep


\section{Examples}\label{examples} 

We give several examples and non-examples of MDS (along with all the possible different ways in which the MDS property can fail). In Example \ref{rnc} we show how the property of being a MDS is neither an open, nor a closed condition.


\begin{example}
Projective $\QQ$-factorial toric varieties are MDS, as they have Cox rings which are polynomial algebras generated by sections corresponding to the 
$1$-dimensional rays of the defining fan \cite{Cox}. 
\end{example}

\begin{example}
$\QQ$-factorial varieties of Fano type are MDS if $\ch \kk=0$ \cite{BCHM}. A variety $X$ is said to be \emph{of Fano type} if there is a Kawamata log-terminal (klt) pair $(X,\De)$, such that $-(K_X+\Delta)$ is ample. Examples include toric varieties, Fano varieties ($\De=\emptyset$) and weak Fano varieties ($-K_X$ is big and nef) with klt singularities. SQMs of varieties of Fano type are of Fano type in characteristic zero (see for example \cite{GOST}, \cite{KO}). 
\end{example}


\begin{example}
Any projective $\QQ$-factorial variety with $\rho=1$ is trivially a MDS. Starting with $\rho\geq2$, there is no classification for MDS, not even for rational surfaces (see Sections \ref{section wpp} and \ref{toricPicard1}).
\end{example}

\begin{example} 
A projective, normal, $\QQ$-factorial surface $X$ is a MDS if and only if the Mori cone $\NE(X)$ is rational polyhedral and every nef divisor $D$ is semiample.
By Zariski's theorem \cite[Rmk. 2.1.32]{Laz}, every movable divisor on a projective surface is semiample. In particular, $\Mov(X)=\Nef(X)$. Hence, a nef divisor $D$ is semiample if and only if a multiple $mD$ is movable for some $m>0$.  
\end{example}


\begin{example}\label{Hesse}
Let $X$ be the blow-up of $\PP^2$ at points $p_1,\ldots, p_r$  in general position. If $r\leq8$, $X$ is a del Pezzo surface $\NE(X)$
is generated by the (finitely many) $(-1)$-curves if $r\geq3$. It follows by induction on $r$ that every nef divisor is semiample.

If $r\geq9$ and the points $p_1,\ldots, p_r$  are in very general position, then $X$ has infinitely many $(-1)$-curves (hence, $\Eff(X)$ has infinitely many extremal rays and $X$ is not a MDS). It is enough to prove that there are  infinitely many $(-1)$-classes when $r=9$ and the points are the base points of a general cubic pencil.
In this case
$$\phi_{|-K_X|}:X\ra\PP^1$$
is an elliptic fibration whose sections are the $(-1)$-curves on $X$. Sections of $\phi$ correspond to $k(t)$-points of the generic fiber $E=X_{k(t)}$ (an elliptic curve over $k(t)$). The Mordell-Weil group $\Pic^0(E)$ is the group of sections of $\pi$, once we fix one section as the identity. It follows that $\Pic^0(E)$ is infinite if for a smooth cubic $C$ containing $p_1,\ldots, p_9$ if $\cO(p_i-p_j)\in\Pic^0(C)$ is non-torsion line bundle for some $i\neq j$. 

When $X$ contains only finitely many $(-1)$- curves (an \emph{extremal} rational elliptic surface), $X$ is a MDS \cite{AL_anticanonical}.
There is a complete classification extremal rational elliptic surfaces, by Miranda-Persson in characteristic zero \cite{MP} and Lang in positive characteristic \cite{Lang1, Lang2}. For example, it follows from this classification that if $\ch \kk\neq 2,3,5$ then the blow-up $X$ of $\PP^2$ at distinct points  $p_1,\ldots, p_9$ which are the base points of a cubic pencil, is extremal if and only if the points are the $9$ flexes of a  smooth cubic in the pencil, i.e., this is the Hesse configuration in $\PP^2$ (unique, up to $\PGL_3$). 
\end{example}

\

\begin{example}\label{Hilbert}
Let $X$ be the blow-up of $\PP^n$ at very general points $p_1,\ldots, p_r$ and let $E_1,\ldots, E_r$ be the corresponding exceptional divisors. Generalizing the case of del Pezzo surfaces, the following are equivalent  \cite{Mukai}, \cite{CT1}:

\bi
\item[(a) ] $X$ is a MDS
\item[(b) ] $\Eff(X)$ is rational polyhedral\footnote{$\Nef(X)$ is rational polyhedral, generated by semiample divisors for $r\leq 2n$.}; 
\item[(c) ] The following inequality holds:
$$\frac{1}{n+1}+\frac{1}{r-n-1}>\frac{1}{2}.$$
\ei

The Weyl group $\cW$ associated to the three-legged Dynkin diagram 
$T_{2,n+1,r-n-1}$ acts on $\Pic(X)$ preserving effective divisors. Every element in the orbit $\cW.E_1$ (which contains all $E_i$'s) generates and
extremal ray of $\Eff(X)$. The group $\cW$ is finite if and only if the above inequality holds, which for $n\geq5$ translates to $r\leq n+3$.

Assume $r=n+3$. Let $C$ be the unique rational normal curve in $\PP^n$ passing through $p_1,\ldots,p_{n+3}$. Then $X$ is a moduli space of parabolic rank $2$  vector bundles on $(C, p_1,\ldots,p_{n+3})$  \cite{Bauer}, \cite{MukaiBook}, \cite{Mukai}.  Varying stability gives rise to all the SQMs of $X$. In particular, $X$ has an SQM which is a weak Fano, hence, $X$ is of Fano type (see also \cite{AM}). 
\end{example}

\begin{example}\label{rnc}
Generalizing Ex. \ref{Hilbert} for $r=n+3$, let $X$ be the blow-up of $\PP^n$ at any number $r$ of points lying on on a rational normal curve. 
Then $X$ is a MDS \cite{CT1}. Hence, being a MDS is not an open condition. We now give an example (due to Hassett and Tschinkel) that shows that being a MDS is not a closed condition either. 

Consider a family of blow-ups $X_t$ of $\PP^3$ along points $p^t_1,\ldots, p^t_9$ lying on some rational normal curve (hence, $X_t$ is a MDS). Such a family admits a degeneration to the blow-up $X_0$ of $\PP^3$ at nine points which are the intersection points of two smooth cubics contained in a plane $\La\subset\PP^3$ (we may assume that the nine points are not the nine flexes of the cubics). Let  $E_1,\ldots,E_9$ be the exceptional divisors on $X_0$ and let $S$ be the proper transform of the plane. As $X_0$ is an equivariant $\GG_a$-compactification of $\PP^3\setminus\La=\GG_a^3$, it follows that  $\Eff(X_0)$ is generated by 
$E_1,\ldots,E_9$, while  $\NE(X_0)$ is generated by curves in $S$. As the restriction map $\Pic(S)\ra\Pic(X_0)$ is an isomorphism, it follows that 
$\NE(X_0)=\NE(S)$ via this identification. As seen in Ex. \ref{Hesse}, $\NE(S)$ is not a rational polyhedral cone if the cubic pencil is not the Hesse pencil. Hence, $X_0$ is not a MDS. 
\end{example}

\begin{example}
If $X$ is a Calabi-Yau variety of dimension at most $3$, then $X$ is a MDS if and only if $\Eff(X)$ is rational polyhedral, generated by effective divisor classes. (The abundance conjecture implies the same statement in higher dimensions \cite[Cor. 4.5]{McK_survey}.) This is clearly the case if $\rho(X)=1$.
If $X$ is a K3 surface with $\rho(X)\geq3$, $\Eff(X)$ is rational polyhedral if and only if $\Aut(X)$ is finite 
(\cite[Thm. 1, Rmk. 7.2]{Kovacs}, \cite{PS-S}). In this case, $\Eff(X)$ is generated by smooth rational curves. 
If $\rho(X)=2$, although  $\Eff(X)$ is rational polyhedral, it may not be generated by effective classes \cite[Thm. 2]{Kovacs}.
\end{example}

\begin{example}\label{toric_generalizations}
Rational normal projective varieties with  a \emph{complexity one torus action} are MDS by \cite{HaussenSuss}. 
Such varieties $X$ admit a faithful action of a torus of dimension $\dim(X)-1$. Examples include projectivizations of toric  rank $2$ vector bundles (see \ref{toricvb}) and several singular del Pezzo surfaces. 

By \cite{Brion}, \emph{wonderful varieties} are MDS. Wonderful varieties admit an action of a semi-simple algebraic group $G$ which has finitely many orbits.
Examples include toric varieties, flag varieties $G/P$ and  and the complete symmetric varieties of De Concini and Procesi \cite{DeConciniProcesi_symmetric}.
\end{example}


\section{Structure Theory}\label{structure} 

As for log-Fano varieties, there is little ``structure theory" for MDS:
\begin{itemize}
\item If $X$ is a MDS, any normal projective variety which is an SQM of $X$, is also a MDS. This follows from the fact that the $f_i$ of Def. 
\ref{mds} are the only SQMs of $X$ (see Rmk. \ref{rmks1}).

\item \cite[Thm. 1.1]{Okawa} If $f: X\ra Y$ is a  surjective morphism of projective normal $\QQ$-factorial varieties and $X$ is a MDS, then $Y$ is a MDS. 
When $f$ is birational, this follows from \cite{HK} (see Rmk. \ref{rmk2}).
\end{itemize}

\subsection{\bf Projective bundles.}
The projectivization $\PP(E)$ of a vector bundle $E$ on a MDS may or may not be a MDS. 
\subsubsection{}
If $L_1,\ldots, L_k$ are line bundles on a MDS $X$, then $\PP(L_1\oplus\ldots\oplus L_k)$ is also a MDS \cite[Thm. 3.2]{Brown}, \cite[Prop. 2.6]{CG}
(see also \cite{Jow}).  

\subsubsection{Toric vector bundles.}\label{toricvb}
A vector bundle $E$ on a toric variety $X$ is called toric if $E$ admits an action of the open torus of $X$ that is linear on fibers and compatible with the action on the base. For example, a direct sum of line bundles is a toric vector bundle. By \cite{GHPS}, a projectivized toric bundles $\PP(E)$ is a MDS if and only if a certain blow-up $Y$ of the fiber of $\PP(E)\ra X$ above the identity point of the torus is a MDS.  Hence, toric $\PP^1$-bundles are always MDS  (see also Ex. \ref{toric_generalizations}). In fact, any blow-up of a projective space along linear subspaces can appear as the variety $Y$ \cite[Cor. 3.8]{GHPS}
(in particular, Ex. \ref{Hilbert}, Ex. \ref{rnc}). Moreover, there is  an example of a toric vector bundle on the Losev-Manin space $\LM_n$ such that $Y =\M_{0,n}$ \cite[p. 21]{GHPS}  (see \ref{LM} for details on Losev-Manin spaces). 

\

The question whether $\PP(E)$ is a MDS seems difficult for non-toric vector bundles $E$, even when $\rk E=2$  \cite{MOS}. 

\subsection{Ample divisors} An ample divisor in a MDS may or may not be a MDS. A question of Okawa: does every MDS have a (not necessarily ample) divisor which is a MDS?

\subsubsection{\bf Lefschetz-type theorems \cite{Jow}}
If $X$ is a smooth MDS of dimension $\geq4$ over $\CC$ which satisfies a certain GIT condition, 
then any smooth ample divisor $Y\subset X$ is a MDS. Moreover, the restriction map identifies $\NN^1(X)$ and $\NN^1(Y)$. Under this identification,  every Mori chamber of $Y$ is a union of some Mori chambers of $X$ and $\Nef(Y)=\Nef(X)$. The GIT condition is stable under taking products and taking the projective bundle of the direct sum of at least three line bundles. The GIT condition is satisfied by smooth varieties of dimension at least $2$ and with $\rho=1$. For  toric varieties, the GIT condition is equivalent to the corresponding fan $\Si$ being \emph{2-neighborly}, i.e., for any $2$ rays of 
$\Si$, the convex cone spanned by them is also in $\Si$. See also \cite{AL_hypers} for examples of non-ample divisors which are MDS. 

\subsubsection{\bf Hypersurfaces in $\PP^m\times \PP^n$ \cite{Ottem}} 
If $X\subset\PP^n\times \PP^m$ is a hypersurface of type $(d,e)$, the cones $\Nef(X)$, $\Mov(X)$ and $\Eff(X)$ are rational polyhedral. 
If $m,n\geq2$, $X$ is a MDS (as proved also in \cite{Jow}). If $m=1$ and $d\leq n$ or $e=1$, then $X$ is a MDS. However, a very general hypersurface $X\subset \PP^1\times\PP^n$ of degree $(d,e)$ with $d\geq n+1$ and $e\geq2$ is not MDS, as $\Nef(X)$ is generated by $H_1$ and $neH_2-dH_1$ (where $H_i=p_i^*\cO(1)$ and $p_1, p_2$ are the two projections), and the divisor $neH_2-dH_1$ has no effective multiple. 
 As noted in \cite{Ottem}, it is the value of $d$, rather than $-K_X$, that determines whether a general hypersurface of degree $(d, e)$ is a MDS or not. In particular, it is not true that a sufficiently ample hypersurface in a MDS is again a MDS.

\subsection{Smooth rational surfaces}\label{surfaces} 
A smooth rational surface $X$ whose anticanonical class $-K_X$ is big (the Iitaka dimension $\kappa(-K_X)$ is $2$) is a MDS \cite[Thm. 1]{TVV}\footnote{There is evidence that the same result holds for all projective $\QQ$-factorial rational surfaces - see Thm. \ref{-K big}.}. There are examples of smooth rational surfaces with $-K_X$ big, which are not of Fano type \cite{TVV}.
Smooth rational surfaces $X$ with  $\kappa(-K_X)=1$ are MDS if and only if $\Eff(X)$ is rational polyhedral \cite{AL_anticanonical}. It is not clear what this condition means in practice.  By Ex. \ref{Hesse}, if $X=\Bl\PP^2_{p_1,\ldots,p_9}$, where $p_1,\ldots,p_9$ are the base points of a cubic pencil, then  
 $X$ is a MDS if and only if $p_1,\ldots,p_9$ are the $9$ inflection points of the cubics in the pencil (the configuration is unique up to $\Aut(\PP^2)$).
 When the points are not the base points of a cubic pencil, it is not clear what the precise condition should be for $X$ to be a MDS. 

When $\kappa(-K_X)\leq 0$, the question is less settled. There exist smooth rational surfaces (of arbitrarily large Picard number) with $\kappa(-K_X)=-\infty$ which are MDS \cite{HwangPark}.


\subsection{Surfaces with $\rho(X)=2$}
The classification of singular rational MDS surfaces with $\rho(X)=2$  is far from settled (see Sections \ref{section wpp} and \ref{toricPicard1}). 
In general, understanding when the blow-up $\Bl_p X$ of a surface $X$ with $\rho(X)=1$ at a general point $p$ is a MDS, is related to the rationality of Seshadri constants (see Section \ref{Picard1}) and is not understood in most cases. 


\subsection{\bf Singularities of Cox rings and positivity of $-K_X$} 
Assume $\ch \kk=0$ and let $X$ be a MDS. Then $X$ is of Fano type (resp., Calabi-Yau type) if and only if $\Spec(\Cox(X))$ has klt singularities (resp. log canonical singularities) \cite{KO} (see also  \cite{GOST}, \cite{Brown}).  Recall that $X$ is said to be \emph{of Calabi-Yau type} if there exists a log-canonical pair $(X,\De)$ such that $(K_X+\De)$ is $\QQ$-linearly trivial. It would be interesting if the condition $-K_X\in\Eff(X)$ is also reflected in $\Cox(X)$.



\section{Blow-ups of surfaces of Picard number one}\label{Picard1} 

Let $X$ be a projective, $\QQ$-factorial, normal surface with $\rho(X)=1$. Let $H$ be an ample $\QQ$-divisor on $X$ and let
$$w:=H^2.$$

If $p\in X$ is a general point, let $\Bl_p X$ denote the blow-up of $p$ and $E$ be the exceptional divisor. 
The Mori cone  of $\Bl_p X$ has the form
$$\NE(\Bl_p X)=\RR_{\geq0}\{E, R\},\quad R=H-\epsilon E, \quad \epsilon\in\RR_{>0}.$$ 

There are two possibilities: either $R^2=0$, or $R^2<0$. Assume that $R^2=0$. Then $\epsilon=\sqrt{w}$ and we have
$$\Nef(X)=\RR_{\geq0}\{H, R\}.$$
In particular, $\epsilon$ is the \emph{Seshadri constant} $\epsilon(H,p)$ of $H$ at the point $p$. 
Then $\Bl_p X$ is a MDS if and only if $R$ is semiample (in particular, 
$\epsilon\in\QQ$). There are no known examples (in any dimension) of irrational Seshadri constants at points. For example, if $X\subset\PP^3$ is a general quintic surface, it is expected that $\epsilon(\cO(1), p)=\sqrt{5}$ for a general point $p$. 
We discuss other conjectural examples of irrational Seshadri constants in Section \ref{section wpp}.

Assume now $R^2<0$. Then there exists an irreducible curve $C$ on $\Bl_p X$ such that $C^2<0$ and $C$ spans the same ray as $R$. 
Then  $\Bl_p X$ is a MDS if and only if the class 
$$R^{\perp}:=H-\frac{w}{\epsilon}E$$ is semiample, or equivalently, using Zariski's theorem, the ray spanned by $R^{\perp}$ contains a movable divisor. As $E$ and $C$ span 
$\NE(\Bl_p X)$ and $R^{\perp}$ is the extremal ray of $\Nef(X)$, it follows that $R^{\perp}$ is semiample if and only if $C$ is not contained in the base locus of $d(R^{\perp})$, for some $d>0$. We state this observation as a Lemma:

\begin{lemma}\label{cutkosky}
Let $X$ be a projective, $\QQ$-factorial surface with Picard number $\rho(X)=1$ and let $p\in X$ be a general point. Let $\Bl_p X$ be the blow-up of $X$ at $p$ and let $E$ be the exceptional divisor. Assume that $\Bl_p X$ contains an irreducible curve $C\neq E$ such that $C^2<0$. Then $\Bl_p X$ is a MDS if and only if there exists an effective divisor $D$ on $\Bl_p X$ such that $D\cdot C=0$ and the linear system $|D|$ does not contain $C$ as a fixed component. Equivalently, there exists a curve $\bar D$ on $X$ that intersects the image $\bar C$ of $C$ in $X$ only at $p$ and with multiplicity one. 
\end{lemma}

\begin{rmk}\label{artin}
Assume the situation in Lemma \ref{cutkosky} and $\ch \kk>0$.
If $X$ and $p$ can be defined over the algebraic closure of a finite field, then a divisor $D$ as in the Lemma always exists. This follows from \cite{Artin} if $X$ is smooth. In general, one can consider the desingularization of $X$ and the same conclusion holds. 
\end{rmk}


\section{Blow-ups of weighted projective planes}\label{section wpp}
Let $a, b, c>0$ be pairwise coprime integers and consider the weighted projective space
$$\PP=\PP(a,b,c)=\Proj S,$$
where $S=\kk[x,y,z]$ and $x,y,z$ have degrees
$$\deg(x)=a,\quad \deg(y)=b,\quad \deg(z)=c.$$ 

Then $\PP$ is a  toric, projective,  $\QQ$-factorial surface with Picard number one. 
Note that $\PP$ is smooth outside the three torus invariant points, but singular at some of these points if $(a,b,c)\neq(1,1,1)$. 
If $D_1, D_2, D_3$ are the  torus invariant (Weil) divisors, let  
$$H=m_1 D_1+m_2 D_2+m_3 D_3,$$
for some integers $m_1, m_2, m_3$ such that $m_1a+m_2b+m_3c=1$. Then  
$$\Cl(\PP)=\ZZ\{H\},\quad \Pic(\PP)=\ZZ\{abc H\},$$
$$H^2=\frac{1}{abc}.$$ 
Moreover, $\cO_{\Proj S}(d)\cong\cO(dH)$ for all $d\in\ZZ$ and $\HH^0(\PP,\cO(d))$ can be identified with the degree $d$ part $S_d$ of $S$. 
If $\pi: \Bl_e\PP\ra\PP$ is the blow-up map, let  $E=\pi^{-1}(e)$. We abuse notations and denote by $H$ the pull-back $\pi^{-1}(H)$ (note that $e$ does not belong to the support of $H$). We have $\Cl(\Bl_e\PP)=\ZZ\{H, E\}$ and hence  a Cox ring of $\Bl_e\PP$ is
$$\Cox(\Bl_e\PP)=\bigoplus_{d,l\in\ZZ}\HH^0(X, \cO(dH-lE)). $$

It was observed by Cutkosky \cite{Cutkosky} that finite generation of $\Cox(\Bl_e\PP)$ is equivalent to the finite generation of the 
symbolic Rees algebra $R_s(\gp)$ of the prime ideal $\gp$ of $S$ defining the point $e$, or equivalently $\gp$ is a \emph{monomial prime}, i.e.,
the kernel of the $k$-algebra homomorphism:
$$\phi: k[x, y, z]\rightarrow k[t], \quad\phi(x=t^ a,\quad\phi(y)=t^b,\quad\phi(z)=t^c.$$

The  \emph{symbolic Rees algebra} of a prime ideal $\gp$ in a ring $R$, is the ring
$$R_s(\gp):=\bigoplus_{l\geq0}\gp^{(l)},\quad \text{ where }\quad \gp^{(l)}=\gp^lR_{\gp}\cap R.$$ 

In our situation, symbolic Rees algebra $R_s(\gp)$ can be identified with the following subalgebra of $\Cox(X)$:
$$\bigoplus_{d,l\in\ZZ_{\geq0}}\HH^0(X, \cO(dH-lE)),$$
which is clearly finitely generated if and only if  $\Cox(\Bl_e\PP)$ is finitely generated (or equivalently Noetherian).
 
The study of the symbolic Rees algebras $R_s(\gp)$ for monomial primes has a long history: \cite{Huneke1}, \cite{Huneke2}, \cite{Cutkosky}, \cite{GNS}, \cite{GNS2}, \cite{Srinivasan}, \cite{GM}, \cite{GNW},  \cite{KuranoMatsuoka}, \cite{CutkoskyKurano}, \cite{GK}.
Prior to \cite{GK}, the only non-finitely generated examples known were the following: 
\begin{thm}\cite[Cor. 1.2, Rmk. 4.5]{GNW}\label{GNW} Assume $(a,b,c)$ is one of the following:
\bi
\item $(7m-3, 5m^2-2m, 8m-3)$,  with $m\geq 4$ and $3\nmid m$, 
\item $(7m-10, 5m^2-7m+1, 8m-3)$,  with $m\geq 5$, $3\nmid 7m-10$ and $m\not\equiv -7(\textrm{mod } 59)$.
\ei
Then $\Bl_e\PP(a,b,c)$ is not a MDS when $\ch \kk=0$.
\end{thm}

The original proof of Theorem \ref{GNW} involved a reduction to positive characteristic. Using methods of toric geometry, Gonzalez and Karu \cite{GK} gave a different proof to Theorem \ref{GNW}, which allows allows for many more examples of toric surfaces $X$ with Picard number one for which $\Bl_e X$ is not a MDS in characteristic zero (Thm. \ref{GK} - to be discussed in detail in  Section \ref{toricPicard1}). In particular:

\begin{thm}\cite{GK}\label{GK particular}
If $\ch \kk=0$ $\Bl_e\PP(a,b,c)$ is not a MDS if $(a,b,c)$ is one of the following: 
$$(7, 15, 26),\quad (7,17, 22),\quad (10,13, 21),\quad (11, 13, 19),\quad (12, 13, 17).$$
\end{thm}
The above are all the triples  $(a,b,c)$  with $a+b+c\leq50$ that satisfy the conditions in Thm. \ref{GK}. 
Key in all the examples in \cite{GK} is that $\Bl_e \PP$ has a negative curve, other than $E$ (hence, Lemma \ref{cutkosky} applies).

\begin{ques}\label{irrational}
Are there any triples $(a,b,c)$ for which $\sqrt{abc}\notin\ZZ$ and $\Bl_e\PP(a,b,c)$ contains no curves $C\neq E$ with $C^2<0$?
\end{ques}

As explained in Section \ref{Picard1},  if $\sqrt{abc}\notin\ZZ$ and $\Bl_e\PP(a,b,c)$ has no negative curves, then  $\Bl_e\PP$ is not a MDS (in any characteristic), as $\NE(\Bl_e \PP)$ and $\Nef(\Bl_e \PP)$ have an irrational extremal ray generated by $H-\frac{1}{\sqrt{abc}}E$. In particular, Seshadri constant $\epsilon(H,e)$ is irrational. Furthermore, if $\kk=\CC$, the Nagata conjecture for $\PP^2$ and $abc$ points holds \cite[Prop. 5.2.]{CutkoskyKurano}.

\medskip

If $\ch \kk>0$ and $\Bl_e\PP$ is not a MDS, then $\Bl_e\PP$ has no negative curve, 
other than $E$ (see Rmk. \ref{artin}). In particular, either  $\sqrt{abc}\notin\ZZ$ or $H-\frac{1}{\sqrt{abc}}E$ is not semiample.  
If  $\Bl_e\PP(a,b,c)$ has no negative curve in characteristic $p$, by standard reduction $p$ methods, it follows  $\Bl_e\PP(a,b,c)$ has no negative curves in characteristic zero. 

\medskip

\begin{ques}\label{9,10,13}\cite{KuranoMatsuoka}
Does $\Bl_e\PP(9,10,13)$ contain a curve $C\neq E$ with $C^2<0$?
\end{ques}

In Section \ref{toricPicard1} we discuss an approach (for $\ch \kk=0$) to the classifcation problem \ref{classify} by reducing the question to an interpolation problem. In particular, Question \ref{9,10,13} has a negative answer (in $\ch \kk=0$, hence, also in $\ch \kk=p$ for all but finitely many primes $p$) if and 
only if there is an affirmative answer to the following:

\begin{ques}(Question \ref{interpolation_9,10,13})
Let $\De$ be the polytope in $\RR^2$ with coordinates $(0,0)$, $(10,40)$, $(36,27)$. For every $q\geq 1$, let 
$$m_q=\lfloor q\sqrt{1170}\rfloor+1.$$
Is it true that for every $q\geq1$ and any point $(i,j)\in q\De\cap\ZZ^2$, there exists a curve $C\subset\RR^2$ of degree $m_q$ passing through all the points 
$(i',j')\neq(i,j)$ in $q\De\cap\ZZ^2$, but not $(i,j)$?
\end{ques}

Computer calculations show that the answer is affirmative for $q\leq5$. 

\

Most known affirmative results are covered by the following:
\begin{thm}\cite{Cutkosky}\label{-K big}
If the anticanonical divisor of $\Bl_e\PP(a,b,c)$ 
$$-K=(a+b+c)H-E$$ 
is big, then $\Bl_e\PP(a,b,c)$ is a MDS. In particular, if $(-K)^2>0$, i.e., if
$$a+b+c>\sqrt{abc},$$
then $\Bl_e\PP(a,b,c)$ is a MDS.
\end{thm}

Note that if $(a,b,c)\neq(1,1,1)$ and $-K$ is big,  $\Bl_e\PP(a,b,c)$ has a negative curve, other than $E$. 
Several particular cases of Thm. \ref{-K big} were proved previously by algebraic methods \cite{Huneke1},  \cite{Huneke2}. Srinivasan \cite{Srinivasan} gave examples of  triples $(a,b,c)$ for which $\Bl_e\PP(a,b,c)$ is a MDS, but $-K$ is not always big:
\bi
\item[(a) ] $(6,b,c)$, for any $b, c$
\item[(b) ] $(5,77, 101)$ (in this case $\kappa(-K)=-\infty$).
\ei

A particular case of Theorem \ref{-K big} is when one of $a,b,c$ is $\leq4$. As noted in \cite{Cutkosky}, when compared with (b) above, this raises the question whether $\Bl_e\PP(5,b,c)$ is always a MDS.

\section{Blow-ups of higher dimensional toric varieties}\label{higher}

Recall that a toric variety $X$ corresponds to the data $(N,\Si)$ where $N$ is a lattice (a finitely generated free $\ZZ$-module) and a fan $\Si\subset N_{\RR}$. Then $X=X(N,\Si)$ is $\QQ$-factorial if and only the fan $\Si$ is simplicial. Two toric varieties $X=X(N,\Si)$ and $X'=X(N',\Si')$ are isomorphic in codimension one if and only if $\Si$ and $\Si'$ have the same rays. To reduce dimensions when considering Question \ref{Jenia}, one has the following result:

\begin{prop}\cite[Prop. 3.1]{CT4}\label{toric}
Let $\pi:\,N\to N'$ be a surjective map of lattices 
with kernel of rank $1$ spanned by a vector $v_0\in N$. 
Let $\Ga$ be a finite set of rays in $N_{\RR}$ spanned by elements of $N$, 
such that the rays $\pm{R_0}$ spanned by $\pm{v_0}$ are not in~$\Gamma$. Let $\Si'\subset N'_\RR$ be a complete simplicial fan with rays given by 
$\pi(\Ga)$. 
Suppose that the corresponding toric variety $X'$ is projective. Then

\bi
\item[(1) ] There exists a complete simplicial fan $\Si\subset N_\RR$
with rays given by $\Ga\cup\{\pm R_0\}$ and such that the corresponding toric variety $X$ is projective and $\pi$ induces a surjective morphism
 $p:\,X\ra X'$. 

\item[(2) ]  There exists an SQM $Z$ of $\Bl_eX$ such that the rational map $Z\dashrightarrow\Bl_eX'$ induced by $p$ is regular.  In particular, if $\Bl_eX$ is a MDS then $\Bl_eX'$ is a MDS.
\ei
\end{prop}


\begin{cor}\label{project}
Assume $X=X(N,\Si)$ is a toric variety of dimension $n$.  Assume there 
exists a saturated sublattice $$N'\subset N,\quad \rk N'=n-2$$ with the following properties:
\bi
\item[(1) ] The vector space $N'\otimes\QQ$ is generated by rays $R$ of $\Si$ with the property that $-R$ is also a ray of $\Si$.
\item[(2) ] There exist three rays of $\Si$ with primitive generators $u, v, w$ whose images generate
$N/N'$ and such that $$au + bv + cw = 0 \quad (\text{mod } N')$$ for some pairwise coprime integers $a,b,c > 0$. 
\ei
Then there exists a rational map $\Bl_e X\dra\Bl_e\PP(a,b,c)$ which is a composition of SQMs and surjective morphisms between normal, projective, $\QQ$-factorial varieties. In particular, if  $\Bl_e X$ is a MDS, then $\Bl_e\PP(a,b,c)$ is a MDS. 
\end{cor}

\say{\bf Losev-Manin spaces.}\label{LM}
Let $\LM_n$ be the Losev-Manin space \cite{LM}. The space $\LM_n$ can be described also as the blow-up of $\PP^{n-3}$ at points $p_1\ldots, p_{n-2}$ in linearly general position and the proper transforms of all the linear subspaces spanned by the points, in order of increasing dimension. 
The space $\LM_n$ is a toric variety and its fan $\Si$ is the barycentric subdivision of the fan of $\PP^{n-3}$. It has lattice
$$N=\ZZ\{e_1,\ldots,e_{n-2}\}/\ZZ\{e_1+\ldots+e_{n-2}\},$$
and rays generated by the primitive lattice vectors 
$$\sum_{i\in I}e_i,\quad \text{for all } I\subset \{1,\ldots, n-2\}, \quad 1\le \#I\le n-3.$$ 

Notice that rays of this fan come in opposite pairs. To construct, for all $n$, a sublattice $N'\subset N$ satisfying the conditions in Cor. \ref{project}, we can proceed as follows: we partition 
$$ \{1,\ldots, n-2\}=S_1\coprod S_2\coprod S_3$$ 
into subsets of size $a+2, b+2, c+2$ (so $n=a+b+c+8$).
We also fix some indices $n_i\in S_i$, for $i=1,2,3$. 
Let $N'\subset N$ be the sublattice generated by the following vectors:
$$e_{n_i}+e_r\quad\hbox{\rm for}\quad r\in S_i\setminus\{n_i\},\ i=1,2,3.$$
If $\pi: N\ra N/N'$ is the projection map,  then we have the following:
\begin{enumerate}
\item $N'$ is a lattice generated by the vectors $\pi(e_{n_i})$, for $i=1,2,3$;
\item $a\pi(e_{n_1})+b\pi(e_{n_2})+c\pi(e_{n_3})=0$.
\end{enumerate}

\begin{cor}\label{sum}
Let $n=a+b+c+8$, where $a,b,c$ are positive pairwise coprime integers.
If $\Bl_e\LM_n$ is a MDS then $\Bl_e\PP(a,b,c)$ is a MDS. 
\end{cor}

Cor. \ref{sum} and Theorems \ref{GNW} and \ref{GK particular} give examples of integers $n$ when $\Bl_e\LM_n$ is not a MDS (for $n\geq134$ if ones uses 
Theorem \ref{GNW} and $n\geq50$ if one uses Theorem  \ref{GK particular}). A smaller 
 $n$ for which  $\Bl_e\LM_n$ is not a MDS was subsequently obtained in \cite{GK}, using Cor. \ref{project} and projecting from a different sublattice $N'$ used in the proof of Cor. \ref{sum}:
\begin{thm}\cite{GK}\label{n=13}
If $\ch \kk=0$, $\Bl_e\LM_{13}$ is not a MDS. 
\end{thm}
Cor. \ref{project} is used to prove that if $\Bl_e\LM_{13}$ is a MDS, then $\Bl_e\PP(7,15,26)$ is a MDS. 
However,  $\Bl_e\PP(7,15,26)$ is not a MDS (Thm. \ref{GK particular}).
The smallest known $n$ (as of the time of this writing) for which  $\Bl_e\LM_n$ is not a MDS was recently obtained in \cite{HKL} by again using Cor. \ref{project} and projecting from a yet different sublattice:
\begin{add}\cite{HKL}\label{n=10}
If $\ch \kk=0$, $\Bl_e\LM_{10}$ is not a MDS. 
\end{add}

Cor. \ref{project} is used to prove that if $\Bl_e\LM_{10}$ is a MDS, then $\Bl_e\PP(12,13,17)$ is a MDS. 
However,  $\Bl_e\PP(12,13,17)$ is not a MDS (Thm. \ref{GK particular}).

\begin{lemma}
If $\Bl_e\LM_{n+1}$ is a MDS, then $\Bl_e\LM_{n}$ is a MDS.
\end{lemma}

\bp
Note that although there exist forgetful maps $\LM_{n+1}\ra\LM_{n}$, in general it is not clear whether one can resolve the rational map 
$$\Bl_e\LM_{n+1}\dra\Bl_e\LM_{n}$$ by an SQM followed by a surjective morphism. However, if $\Bl_e\LM_{n+1}$ is a MDS, this is always the case, and we are done by \cite{Okawa}.
\ep

As $\Bl_e\LM_{6}$ is a MDS in any characteristic (follows from \cite{C} - see \ref{M0n}; moreover, it is a threefold of Fano type), we are left with:
\begin{ques} 
Is  $\Bl_e\LM_{n}$ a MDS for $7\leq n\leq9$, $\ch \kk=0$?
\end{ques}

\

\say{\bf Losev-Manin spaces and the moduli spaces $\M_{0,n}$}\label{M0n}
There is a close connection between the blow-ups $\Bl_e\LM_{n}$ of the Losev-Manin spaces and the moduli spaces $\M_{0,n}$ of stable, $n$-pointed 
rational curves. By Kapranov \cite{Kapranov}, $\M_{0,n}$ is the blow-up of $\PP^{n-3}$ at points $p_1\ldots, p_{n-1}$ in linearly general position and the proper transforms of all the linear subspaces spanned by the points, in order of increasing dimension. 
Up to changing coordinates, we may assume that 
$$p_1=[1,0,0,\ldots,0],p_2=[0,1,0,\ldots,0],\ldots,p_{n-2}=[0,0,0,\ldots,1],$$
$$p_{n-1}=e=[1,1,1,\ldots,1].$$

Note that $p_{n-1}$ is the identity of the open torus in $\LM_n$. Moreover, $\M_{0,n}$ is the blow-up of $\LM_n$ along $e$, and the (proper transforms of the) linear susbpaces spanned by $e$ and $\{p_i\}_{i\in I}$, for all the subsets $I$ of $\{1,\ldots, n-2\}$ with $1\leq\#I\leq n-5$. In particular, there is a projective birational morphism $\M_{0,n}\ra\Bl_e\LM_{n}$.

\begin{thm}\cite{CT4}\label{main CT}
\bi
\item[(1) ] If $\M_{0,n}$ is a MDS, then $\Bl_e\LM_{n}$ is a MDS;
\item[(2) ]  If  $\Bl_e\LM_{n+1}$ is a MDS, then $\M_{0,n}$ is a MDS. 
\ei
\end{thm}

The existence of forgetful maps $\M_{0,n+1}\ra\M_{0,n}$ implies that if $\M_{0,n+1}$ is a MDS, then $\M_{0,n}$ is a MDS. 
Combined with Cor. \ref{sum} and the resuts in \ref{LM}, Thm. \ref{main CT} gives a negative answer to the question of Hu and Keel \cite{HK} whether
$\M_{0,n}$ is a MDS.
\begin{thm}\cite{CT4},\cite{GK},\cite{HKL}
If $n\geq10$, $\M_{0,n}$ is not a MDS in characteristic $0$.
\end{thm}

Note that $\M_{0,6}$ is a MDS in any characteristic \cite{C} (moreover, it is a threefold of Fano type). The range $7\leq n\leq9$ is still open. 

Part (1) of Thm. \ref{main CT} follows from \cite{HK} (see Rmk. \ref{rmk2}). Part (2) follows from:
\begin{thm}\cite{CT4}\label{Xn}
Let $X_n$ be the toric variety which is the blow-up of $\PP^{n-3}$ along points $p_1,\ldots,p_{n-2}$ and (all the proper transforms of) the linear subspaces of codimension at least $3$ spanned by the points $p_1,\ldots,p_{n-2}$. Then $\Bl_e X_{n+1}$ is an SQM of a $\PP^1$-bundle over $\M_{0,n}$ which is the projectivization of a direct sum of line bundles. 
 \end{thm}

Hence, $\M_{0,n}$ is a MDS if and only if $\Bl_eX_{n+1}$ is a MDS. In particular:
\bi
\item If $n\geq11$, then $\Bl_eX_n$ is not a MDS if $\ch \kk=0$;
\item If $n\leq7$, then $\Bl_eX_n$ is a MDS.
\ei


\say{\bf Further questions.} 
\begin{enumerate}
\item Are there other examples of toric varieties besides Losev-Manin spaces, to which Cor. \ref{project} applies?
\item What are the simplest smooth toric varieties $X$ for which $\Bl_e X$ is not a MDS? Any smooth Fano varieties?
\end{enumerate}

If $X$ is a projective, $\QQ$-factorial toric variety such that all the torus invariant divisors are not movable, then $\Bl_e X$ is not toric. 
It may or may not be a MDS (for example, when $X$ is $\LM_6$ or $\LM_n$ with $n\geq10$). 
If some of the torus invariant divisors are movable, then $\Bl_e X$ may be toric (for example when $X=\PP^n$), but may not even be a MDS (for example, when $X=X_n$ from Thm. \ref{Xn}). It would be interesting to find a geometric criterion for $\Bl_e X$ to not be a MDS.

\section{Blow-ups of toric surfaces}\label{toricPicard1}

In this section we assume 
$$\ch\kk=0.$$
Let $(X_{\De}, H)$ be a polarized  toric projective surface with $H$ an ample $\QQ$-Cartier divisor on $X_{\De}$ corresponding to the rational polytope $\De\subset N^*_{\RR}=\RR^2$. If $X_{\De}$ has Picard number $\rho$, then $\De$ is a rational polytope with $\rho+2$ vertices.  If $d>0$ is an integer such that $d\De$ has integer coordinates, then global sections of $\cO_{X_{\De}}(dH)$ can be identified with Laurent polynomials (considered as regular functions on the open torus):
$$f=\sum_{(i,j)\in d\De\cap\ZZ^2}a_{(i,j)}x^iy^j\in \HH^0(X,\cO(dH)).$$

The vertices of $\De$ correspond to the $\rho+2$ torus invariant points of $X$. A section $f$ vanishes at a torus invariant point if and only if the coefficient $a_{ij}$ of the corresponding vertex in $d\De$ is zero. We fix a vertex $(x_1,y_1)$ of $\De$ and  and let $p_1$ be the corresponding torus invariant point. 
For simplicity, we assume this is the ``leftmost lowest" point of $\De$.

We now translate into linear algebra the condition that a global section of $\cO_{X_{\De}}(dH)$  has a certain multiplicity at the point $e$. Let $N_d$ be the number of lattice points $(i,j)\in d\De\cap\ZZ^2$ and let $R_m$ be the number of derivatives $\de^a_x\de^b_y$ of order $\leq m-1$ in two variables: 
$$R_m=1+2+\ldots+m=\frac{m(m+1)}{2}.$$

\begin{defn}
We order the pairs $(i,j)$ and the pairs $(a,b)$ lexicographically (so the first $(i,j)$ corresponds to the leftmost point $(dx_1,dy_1)$ of $d\De$). We define two $N_d\times R_m$  matrices $A=A_{d,m}$ and $B=B_{d,m}$, whose entries for the pairs $(i,j)$ and $(a,b)$ as are given as follows: 
$$A_{(i,j),(a,b)}=\de^a_x\de^b_y(x^iy^j)(1,1)=a!{i\choose{a}}b!{j\choose{b}},$$
$$B_{(i,j),(a,b)}=i^aj^b.$$
where 
we denote for any integers $n,k$ ($k\geq0$, but $n$ possibly negative) 
$$ {n\choose{k}}=\frac{n(n-1)(n-2)\ldots(n-k+1)}{k!}.$$ 
\end{defn}

We write $N=N_d$, $R=R_m$, $A=A_{d,m}$, $B=B_{d,m}$, whenever there is no risk of confusion.
\begin{lemma}\label{AB}
The matrix $B_{d,m}$ can be obtained from $A_{d,m}$ by a sequence of reversible column operations.
\end{lemma}

\bp 
We claim that for every column $(a,b)$ of $A$, starting from left to right, we can do (reversible) column operations on $A$ involving only previous columns, and end up with the column that has entries $i^aj^b$ for every row $(i,j)$. For simplicity, we may first ignore the $j$'s and consider the situation when one matrix has entries $a!{i\choose{a}}$ and the other $i^a$ (with rows indexed by $i$ and columns by $a$). It is easy to see that one can do reversible column operations from one matrix to the other: use induction on $a$  and  expand the product 
$$i(i-1)(i-2)\ldots(i-a+1).$$
The general case is similar. 

\ep

\begin{lemma}\label{linear algebra}
Let $A$ be an $N\times R$ matrix with entries in $\QQ$. The following are equivalent:
\bi
\item[(a) ] Any linear combination $\sum \al_j R_j$  of the rows $R_j$ of the matrix $A$ that is zero, must have $\al_i=0$,
\item[(b) ] There exists a linear combination of the columns of $A$ that equals the vector $e_i=(0,0,\ldots,1,0\ldots 0)\in\RR^{N_d}$. 
\ei

In particular,  $A$ has rank $N$ if and only if for every $1\leq i\leq N$, there exists a linear combination of the columns of $A$ which equals 
$e_i$. 
\end{lemma}

\bp 
May assume $i=1$. Consider the pairing 
$(,): V\times W\ra\QQ$ with $V$ a $\QQ$-vector space with basis $e_1,\ldots, e_N$ and  $W$ a $\QQ$-vector space with basis $f_1,\ldots, f_R$ and 
$(e_u,f_v)=a_{uv}$. Let $$\phi: V\ra W^*,\quad \phi^*:W\ra V^*$$ be the induced linear maps. Condition (b) is equivalent to the dual vector $e_1^*\in V^*$ being in the image of the map $\phi^*$. Condition (a) is equivalent to  the kernel $K$ of the map $\phi$ being contained in the span of $e_2,\ldots, e_N$. Let 
$I=Im(\phi)\subset W^*$. Hence, there is an exact sequence $$0\ra K\ra V\ra I\ra 0.$$ 
Consider the inclusion map $u: K\ra V$.  Dualizing, it follows that $Im(\phi^*)=I^*=ker(u^*)$. Hence, $e_1^*\in Im(\phi^*)$ if and only if $u^*(e_1^*)=0$. As  $u^*(e_1^*)$ is the linear functional $K\ra \QQ$ given by $k\mapsto e_1^*(k)$, for $k\in K$, it follows that  $u^*(e_1^*)=0$ if and only if $e_1^*(k)=0$, for all $k\in K$, or equivalently, $K$ is contained in the span of $e_2,\ldots, e_N$.
\ep

\begin{lemma} Let $\Bl_e X_{\De}$ be the blow-up of $X_{\De}$ at the identity point $e$ and let $E$ denote the exceptional divisor. 
The following are equivalent:
\bi
\item[(i) ] The linear system  $|dH-mE|$ is empty,
\item[(ii) ] The matrix $A_{d,m}$ has linearly independent rows,
\item[(iii) ] The matrix $B_{d,m}$ has linearly independent rows,
\item[(iv) ] For every $(i,j)\in d\De\cap\ZZ^2$, there exists a polynomial $f(x,y)\in\QQ[x,y]$ of degree $\leq m-1$, such that $ f(i,j)\neq0$ and 
$$f(i',j')=0 \text{ for all } (i',j')\in d\De\cap\ZZ^2, (i'j')\neq(i,j).$$
\ei
\end{lemma}
Equivalently, condition (iv) says that one can separate any lattice point in $d\De$ from the rest by degree $m-1$ plane curves. 

\bp
A non-zero section of $\cO_{X_{\De}}(dH)$ has multiplicity $m$ at the point $e$ if and only if there exists a non-zero linear combination of the rows of  $A_{d,m}$ which is zero. Hence, (i) is clearly equivalent to (ii).  By Lemma \ref{AB}, (ii) is equivalent to (iii). By Lemma \ref{linear algebra}, (iii) is equivalent to (iv).
\ep

\begin{lemma}\label{separate leftmost} 
Let $\Bl_e X_{\De}$ be the blow-up of $X_{\De}$ at the identity point $e$ and let $E$ denote the exceptional divisor. 
The following are equivalent:
\bi
\item[(i) ] All non-zero sections of the linear system $|dH-mE|$ (if any) define curves that pass through the torus invariant point $p_1$, 
\item[(ii) ] There exists a linear combination of the columns of the matrix $A_{d,m}$ that equals the vector $(1,0\ldots, 0)\in\RR^{N_d}$, 
\item[(iii) ] There exists a linear combination of the columns of the matrix $B_{d,m}$ that equals the vector $(1,0\ldots, 0)\in\RR^{N_d}$, 
\item[(iv) ] There exists a polynomial $f(x,y)\in\QQ[x,y]$ of degree $\leq m-1$, such that $f(dx_1,dy_1)\neq0$ and 
$$f(i,j)=0 \text{ for all } (i,j)\in d\De, (i,j)\neq(dx_1,dy_1).$$
\ei
\end{lemma}
Equivalently, condition (iv) says that there exists a plane curve of degree $\leq m-1$ that passes through all the lattice points in $d\De$, except the lefmost point.  

\bp

Condition (i) is equivalent to the fact that any non-zero section of $\cO_{X_{\De}}(dH)$ which has multiplicity $m$ at the point $e$, must have the coefficient 
$a_{(dx_1,dy_1)}$ is zero. Equivalently, any linear combination $\sum \al_i R_i$  of rows $R_i$ of the matrix $A$ that is zero, must have $\al_1=0$. By Lemma \ref{linear algebra} this is equivalent to condition (ii). Lemma \ref{AB} implies that (ii) and (iii) are equivalent. Condition (iv) is just a reformulation of (iii).
\ep

\

Consider now the situation when $\rho(X_{\De})=1$ (i.e.,  $\De$ is a triangle) and $\Bl_e(X_{\De})$ has a curve $C\neq E$ with $C^2<0$. As in \cite{GK}, 
we assume that the point $(0,0)$ is one vertex of $\De$, the point $(0,1)$ lies in the interior of a non-adjacent edge, and moreover, $C$ is the proper transform of the closure $\bar C$ of the curve defined by the section  $1-y$ of $\cO_{X_{\De}}(H)$. Then $\bar C=H$ in $\Cl(X_{\De})$ and 
$$C=H-E\quad \text{in }\Cl(\Bl_eX_{\De}).$$ 
The condition $C^2<0$ is equivalent to $$w:=H^2=2(\text{Area}(\De))<1.$$ 

Denote by $(x_1,y_1)$ the leftmost point of $\De$ and by $(x_2,y_2)$  the rightmost point of $\De$. 
Let $p_1$, respectively $p_2$, be the corresponding torus invariant points.  
Note that $\bar C$ contains $p_1$ and $p_2$. Moreover, $w=H^2=x_2-x_1$ is the width of $\De$.

\medskip 

The main theorem in \cite{GK} becomes an instance of the following more general statement, which shows that the question of $\Bl_e X_{\De}$ not being a MDS is equivalent to solving an interpolation problem for points in the (usual) affine plane.
 
\begin{prop}\label{GK general}
Let $(X_{\De},H)$ be a polarized projective toric surface with $\rho(X_{\De})=1$ corresponding to a triangle $\De$ as above. Assume 
$$w=H^2<1.$$
Then  $\Bl_e X_{\De}$ is not a MDS if and only if for any sufficiently divisible integer $d>0$ such that $d\De$ has integer coordinates, there exists a curve $C\subset\AA^2$ of degree $dw-1$ that passes through all the lattice points $d\De\cap\ZZ^2$ except the point $(dx_1,dy_1)$.
\end{prop}

\bp
By Lemma \ref{cutkosky}, $\Bl_e X_{\De}$ is not a MDS if and only if any non-zero effective divisor $D$ with class $dH-dwE$ ($d>0$)
contains $C$ in its fixed locus, or equivalently,  the image $\bar D$ of $D$ in $X_{\De}$ contains some other point of $\bar C$ than $e$ (for example $p_1$). 
Hence, $\Bl_e X_{\De}$ is not a MDS if and only if for any sufficiently large and divisible $d$, any element of the linear system $|dH-dw E|$ contains $p_1$. The result now follows from Lemma \ref{separate leftmost}. 
\ep

The difficult part is of course to solve the interpolation problem posed in Prop.  \ref{GK general}. We claim that the main theorem in \cite{GK} gives sufficient (but not necessary) conditions for this.

\begin{thm}\cite[Thm. 1.2]{GK}\label{GK}
Let $(X_{\De},H)$ be a polarized projective toric surface with Picard number one, 
corresponding to a triangle $\De$ as above and assume  $$w=H^2<1.$$ 

If $s_1<s_2<s_3$ are the slopes defining the triangle $\De$, let 
$$n=\#([s_1,s_2]\cap\ZZ).$$

Assume that $$\#((n-1)[s_2,s_3]\cap\ZZ)=n,\quad \text{and}\quad ns_2\notin\ZZ.$$ 
Then for any integer $d>0$ such that $d\De$ has integer coordinates,  there exists a curve $C\subset\AA^2$ of degree $dw-1$ that passes through all the lattice points $d\De\cap\ZZ^2$ except the leftmost point $(dx_1,dy_1)$. In particular, $\Bl_e X_{\De}$ is not a MDS by Proposition \ref{GK general}.
\end{thm}

As mentioned in \cite{GK}, $\#([s_1,s_2]\cap\ZZ)$ represents the number of points in $d\De\cap\ZZ^2$ (for any $d$ such that $d\De$ has integer coordinates) lying in the second column from the left, i.e., the column with $x$ coordinate $mx_1+1$. Similarly, for any $k\geq1$, 
the number 
$$\#((k-1)[s_2,s_3]\cap\ZZ)$$ is the number of points in $d\De\cap\ZZ^2$ lying in the $k$-th column from the right, i.e., the column with $x$ coordinate 
$mx_2-(k-1)$. None of these numbers depend on $d$. The condition $ns_2\notin\ZZ$ is equivalent to the $(n+1)$-th column from the right not containing a lattice point on the top edge (see Rmk. \ref{ss3}).

\bp[Proof of Theorem \ref{GK}]
As in \cite{GK}, we first transform the triangle $d\De$ by integral translations and shear transformations $(i,j)\mapsto(i,j+ai)$ for $a\in\ZZ$. Clearly, the assumptions still hold for the new triangle. To see that the conclusion is also not affected, recall that the conclusion is equivalent to the fact that any section $f$ of $\HH^0(X_{\De},dH)$
$$f(x,y)=\sum_{(i,j)\in d\De\cap\ZZ^2} a_{(i,j)}x^iy^j$$ that vanishes to order $dw$ at $e=(1,1)$ has the coefficient  $a_{(dx_1,dx_2)}=0$ (i.e., $f$ vanishes at the torus invariant point $p_1$). The translation operation multiplies $f$ with a monomial, and the shear transformation performs a change of variables on the torus. The two operations do not affect the order of vanishing of $f$ at $e$ or whether $f$ vanishes at $p_1$.

We first apply a shear transformation, so that $-2<s_2<-1$ (possible since $s_2\notin\ZZ$). We then translate the triangle so that the leftmost point moves to a point with $x$-coordinate $-1$ and the rightmost point moves to a point on the $x$-axis. As there are precisely $n$ lattice points in the $n$-th column from the right, it follows from $-2<s_2<-1$ and that the $n$ points are, in coordinates
$$(\al,0),\quad (\al,1),\quad (\al,2),\quad \ldots, (\al,n-1),\quad \text{for some } \al\geq0,$$
along with $0\leq s_3$. Note  $ns_2\notin\ZZ$ implies $\al>0$\footnote{We may also take $\al>0$ at the expense of proving the statement only for sufficiently large and divisible $d$.}.
It also follows that for all $0\leq i\leq n-1$, the column in $d\De$ with $x$-coordinate $\al+i$ has exactly $i$ lattice points:
$$(\al+i,0),\quad (\al+i,1),\quad (\al+i,2),\quad \ldots, (\al+i,n-1-i).$$
We denote these points $\{Q_j\}$ (a total of $\frac{n^2+n}{2}$ points). Let the $n$ lattice points in $d\De$ in the second column from the left   be
$$P_0=(\be,0),\quad P_1=(\be+1,0),,\quad \ldots, P_{n-1}=(\be+n-1,0),$$
for some $\be\geq0$. As $-2<s_2<-1$,  the rightmost point must be
$$L=(-1,\be+n+1).$$
As the width of $d\De$ is $dw$, the integers $\al,\be$ are related to $w, s_2$ by
$$\al=dw-n, \quad \be=-s_2(dw)-n-1,\quad -s_2=\frac{\be+n+1}{\al+n}$$


\begin{lemma}\label{interpolation}
There is a unique curve $C$ of degree $\leq n$ passing through the $\frac{n^2+3n}{2}$ points $\{P_i\}$ and $\{Q_i\}$. The curve $C$ passes through the point $L$ if and only if $n\be=(n+1)\al$ (or, equivalently, $-s_2=1+\frac{1}{n}$). 
\end{lemma}


\begin{rmk}\label{ss3}
It is not hard to see that the condition $ns_2\notin\ZZ$ is equivalent to $-s_2=1+\frac{1}{n}$, which in turn says that 
$(n+1)$-th column from the right not containing a lattice point on the top edge. 
\end{rmk}

Assuming Lemma \ref{interpolation}, Theorem \ref{GK} follows by considering the union $C'$ of the curve $C$ with all the vertical lines 
$$x=1,\quad x=2,\quad\ldots x=(\al-1).$$
Note that the degree of $C'$ equals $dw-1$. Clearly, if $ns_2\notin\ZZ$, Lemma  \ref{interpolation} implies that $C'$ does not pass through $L$. 
 
\bp[Proof of Lemma \ref{interpolation}] 
We first write down a basis $G_0,\ldots, G_n$ for the vector space of polynomials in $\QQ[x,y]$ of degree $\leq n$ that vanish at the points $\{Q_j\}$ as follows. For all  $0\leq i\leq n$, let
$$G_i(x,y)={x-\al\choose{i}}{y\choose{n-i}}.$$

Consider now the equation of a curve $C$ that passes through $\{Q_j\}$:
$$f(x,y)=\sum_{i=0}^n c_{i}G_i(x,y),\quad c_{i}\in\QQ.$$

Let $M$ be the $(n+1)\times (n+1)$ matrix with rows indexed by points $P_0, P_1,\ldots, P_{n-1}, L$ (hence, the last row corresponds to $L$) and columns indexed by  $G_0,\ldots, G_n$, such that the entry corresponding to the row $P_i$ (resp. $L$) and column $G_j$ is $G_j(P_i)$ (resp. $G_j(L)$), i.e., 
$$M_{P_i,G_j}=G_j(0,\be+i)={-\al\choose{j}}{\be+i\choose{n-j}},\quad 0\leq i\leq n-1,$$
$$M_{L,G_j}=G_j(-1,\be+n+1)={-1-\al\choose{j}}{\be+n+1\choose{n-j}}.$$

Let $M'$ be the $n\times (n+1)$ matrix obtained by taking the first $n$ rows of $M$. Clearly, there is a unique curve $C$ passing through $\{Q_j\}$ and $\{P_j\}$ if and only if there is a unique solution ${\bf c}=(c_i)$ (up to scaling) to the linear system $M'\cdot{\bf c}={\bf 0}$, i.e., 
$$\rk M'=n.$$ 

To prove this, successively substract row $P_{n-2}$ from row $P_{n-1}$, row $P_{n-3}$ from row $P_{n-2}$, etc, row $P_0$ from row $P_1$. The result is that the last column of $M'$ has the last $(n-1)$ entries $0$.  Substracting row $P_{n-2}$ from row $P_{n-1}$, row $P_{n-3}$ from row $P_{n-2}$, etc, row $P_1$ from row $P_2$  leaves the second column of $M'$ with the last $(n-2)$ entries $0$. Continuing in the same fashion (and using the relation
${k+1\choose{l+1}}={k\choose{l+1}}+{k\choose{l}}$) we obtain an ``upper diagonal" matrix $n\times (n+1)$ matrix $M''$ with entries:
$$
M''_{P_i,G_j} = \left\{
        \begin{array}{ll}
            {-\al\choose{j}}{\be\choose{n-i-j}} & \quad\text{if} \quad  i+j\leq n \\
            0 & \quad \text{if} \quad  i+j> n
        \end{array}
    \right.
$$

Hence, $\rk M'=\rk M''=n$. 

We now prove that $\det M=0$ if and only if $n\be=(n+1)\al$. Clearly, the curve $C$ passes through the point $L$ if and only if $\det M=0$, hence, this would finish the proof. Let $\tilde{M}$ be the matrix obtained by adding to the matrix $M''$ the last row of $M$, i.e., 
$$
\tilde{M}_{P_i,G_j} = \left\{
        \begin{array}{ll}
            {-\al\choose{j}}{\be\choose{n-i-j}} & \quad\text{if} \quad  i+j\leq n \\
            0 & \quad \text{if} \quad  i+j> n,
        \end{array}
    \right.
$$
$$\tilde{M}_{L,G_j} = {-1-\al\choose{j}}{\be+n+1\choose{n-j}}.$$
Clearly, $\det M=\det\tilde{M}$. Let $\tilde{M}^{(1)}$ be the matrix obtained from $\tilde{M}$ by first dividing the column corresponding to $G_j$ by ${-\al\choose{j}}$ (for every $j$) and multiplying the last row with $\al$. Using that ${-1-\al\choose{j}}={-\al\choose{j}}\frac{\al+j}{\al}$, the entries of $\tilde{M}^{(1)}$ are given by
$$
\tilde{M}^{(1)}_{P_i,G_j} = \left\{
        \begin{array}{ll}
            {\be\choose{n-i-j}} & \quad\text{if} \quad  i+j\leq n \\
            0 & \quad \text{if} \quad  i+j> n,
        \end{array}
    \right.
$$
$$\tilde{M}^{(1)}_{L,G_j} = (\al+j){\be+n+1\choose{n-j}}.$$

Let $\tilde{M}^{(2)}$ be the matrix obtained from $\tilde{M}^{(1)}$ by first multiplying the last row with $(-1)$, then adding to the last row the sum  of rows:
$${n+1\choose{0}}(\text{row } P_0)+{n+1\choose{1}}(\text{row } P_1)+\ldots+{n+1\choose{n}}(\text{row } P_{n-1}),$$
then finally dividing the last row by $(\frac{1}{\be+n+1})$. 
Using the identities
$$\sum_{i=0}^{n-j}{\be\choose{n-i-j}}{n+1\choose{i}}={\be+n+1\choose{n-j}},\quad l{n\choose{l}}=n{n-1\choose{l-1}}$$
it follows that  the entries in the last row of $\tilde{M}^{(2)}$ are:
$$\tilde{M}^{(2)}_{L,G_j} = {\be+n\choose{n-j-1}},\quad 1\leq j\leq n,$$
$$\tilde{M}^{(2)}_{L,G_0} = {\be+n\choose{n-1}}-\frac{(\al+n)(n+1)}{(\be+n+1)}.$$
 
Finally, let $\tilde{M}^{(3)}$ be the matrix obtained from $\tilde{M}^{(2)}$ by
substracting from the last row, the following sum of rows:
$${n\choose{0}}(\text{row } P_1)+{n\choose{1}}(\text{row } P_2)+\ldots+{n\choose{n-2}}(\text{row } P_{n-1}).$$
The matrix $\tilde{M}^{(3)}$ has entries
$$\tilde{M}^{(3)}_{L,G_j} =0,\quad 1\leq j\leq n,$$
$$\tilde{M}^{(3)}_{L,G_0} =n-\frac{(\al+n)(n+1)}{(\be+n+1)}.$$
Note that  $\tilde{M}^{(3)}_{L,G_0}=0$ if and only if $n\be=(n+1)\al$.
As $\tilde{M}^{(3)}$ is an upper triangular matrix with $\det\tilde{M}^{(3)}=\tilde{M}^{(3)}_{L,G_0}$, the result follows. 
\ep
\ep

There are other possible applications of Prop. \ref{GK general} that are not covered by Theorem \ref{GK} towards the classification problem \ref{classify} (see also \cite{He}).  For toric surfaces of higher Picard number, we expect that solving an interpolation problem analogous to the one posed in Prop. \ref{GK general}
will lead to examples of non Mori Dream Spaces. An interesting question is whether there is  higher dimensional version of Prop. \ref{GK general}.



\end{document}